\magnification=\magstep1
\hsize=16.5 true cm 
\vsize=23.6 true cm
\font\bff=cmbx10 scaled \magstep1
\font\bfff=cmbx10 scaled \magstep2
\font\bffg=cmbx10 scaled \magstep3

\font\smc=cmcsc10 
\parindent0cm
\overfullrule=0cm
\def\cl{\centerline}           %
\def\rl{\rightline}            %
\def\bp{\bigskip}              %
\def\mp{\medskip}              %
\def\sp{\smallskip}            %
           %
\def\Bbb#1{\hbox{\boldmas #1}} %
\def\K{\Bbb K}                 %
\def\N{\Bbb N}                 %
\def\Q{\Bbb Q}                 %
\def\R{\Bbb R}                 %
\def\C{\Bbb C}                 %
\expandafter\edef\csname amssym.def\endcsname{%
       \catcode`\noexpand\@=\the\catcode`\@\space}

\catcode`\@=11

\def\undefine#1{\let#1\undefined}
\def\newsymbol#1#2#3#4#5{\let\next@\relax
 \ifnum#2=\@ne\let\next@\msafam@\else
 \ifnum#2=\tw@\let\next@\msbfam@\fi\fi
 \mathchardef#1="#3\next@#4#5}
\def\mathhexbox@#1#2#3{\relax
 \ifmmode\mathpalette{}{\m@th\mathchar"#1#2#3}%
 \else\leavevmode\hbox{$\m@th\mathchar"#1#2#3$}\fi}
\def\hexnumber@#1{\ifcase#1 0\or 1\or 2\or 3\or 4\or 5\or 6\or 7\or 8\or
 9\or A\or B\or C\or D\or E\or F\fi}

\font\tenmsa=msam10
\font\sevenmsa=msam7
\font\fivemsa=msam5
\newfam\msafam
\textfont\msafam=\tenmsa
\scriptfont\msafam=\sevenmsa
\scriptscriptfont\msafam=\fivemsa
\edef\msafam@{\hexnumber@\msafam}
\mathchardef\dabar@"0\msafam@39
\def\dashrightarrow{\mathrel{\dabar@\dabar@\mathchar"0\msafam@4B}}
\def\dashleftarrow{\mathrel{\mathchar"0\msafam@4C\dabar@\dabar@}}

\def\ulcorner{\delimiter"4\msafam@70\msafam@70 }
\def\urcorner{\delimiter"5\msafam@71\msafam@71 }
\def\llcorner{\delimiter"4\msafam@78\msafam@78 }
\def\lrcorner{\delimiter"5\msafam@79\msafam@79 }
\def\yen{{\mathhexbox@\msafam@55}}
\def\checkmark{{\mathhexbox@\msafam@58}}
\def\circledR{{\mathhexbox@\msafam@72}}
\def\maltese{{\mathhexbox@\msafam@7A}}

\font\tenmsb=msbm10
\font\sevenmsb=msbm7
\font\fivemsb=msbm5
\newfam\msbfam
\textfont\msbfam=\tenmsb
\scriptfont\msbfam=\sevenmsb
\scriptscriptfont\msbfam=\fivemsb
\edef\msbfam@{\hexnumber@\msbfam}
\def\Bbb#1{{\fam\msbfam\relax#1}}
\def\widehat#1{\setbox\z@\hbox{$\m@th#1$}%
 \ifdim\wd\z@>\tw@ em\mathaccent"0\msbfam@5B{#1}%
 \else\mathaccent"0362{#1}\fi}

\def\widetilde#1{\setbox\z@\hbox{$\m@th#1$}%
 \ifdim\wd\z@>\tw@ em\mathaccent"0\msbfam@5D{#1}%
 \else\mathaccent"0365{#1}\fi}
\font\teneufm=eufm10
\font\seveneufm=eufm7
\font\fiveeufm=eufm5
\newfam\eufmfam
\textfont\eufmfam=\teneufm
\scriptfont\eufmfam=\seveneufm
\scriptscriptfont\eufmfam=\fiveeufm

\newsymbol\risingdotseq 133A
\newsymbol\fallingdotseq 133B
\newsymbol\complement 107B
\newsymbol\nmid 232D
\newsymbol\rtimes 226F
\newsymbol\thicksim 2373

\font\eightmsb=msbm8   \font\sixmsb=msbm6   \font\fivemsb=msbm5
\font\eighteufm=eufm8  \font\sixeufm=eufm6  \font\fiveeufm=eufm5
\font\eightrm=cmr8     \font\sixrm=cmr6     \font\fiverm=cmr5
\font\eightbf=cmbx8    \font\sixbf=cmbx6    
      \font\eighti=cmmi8   \font\sixi=cmmi6
\font\ninesy=cmsy9     \font\eightsy=cmsy8  \font\sixsy=cmsy6
     \font\eightit=cmti8  
     \font\eightsl=cmsl8  
     \font\eighttt=cmtt8

\font\eightsmc=cmcsc8
\newskip\ttglue
\newfam\smcfam
\def\eightpoint{\def\rm{\fam0\eightrm}%
  \textfont0=\eightrm \scriptfont0=\sixrm \scriptscriptfont0=\fiverm
  \textfont1=\eighti \scriptfont1=\sixi \scriptscriptfont1=\fivei
  \textfont2=\eightsy \scriptfont2=\sixsy \scriptscriptfont2=\fivesy
  \textfont3=\tenex \scriptfont3=\tenex \scriptscriptfont3=\tenex
  \def\smc{\fam\smcfam\eightsmc}
  \textfont\smcfam=\eightsmc          
\textfont\eufmfam=\eighteufm              \scriptfont\eufmfam=\sixeufm
     \scriptscriptfont\eufmfam=\fiveeufm
\textfont\msbfam=\eightmsb            \scriptfont\msbfam=\sixmsb
     \scriptscriptfont\msbfam=\fivemsb
\def\it{\fam\itfam\eightit}%
  \textfont\itfam=\eightit
  \def\sl{\fam\slfam\eightsl}%
  \textfont\slfam=\eightsl
  \def\bf{\fam\bffam\eightbf}%
  \textfont\bffam=\eightbf \scriptfont\bffam=\sixbf
   \scriptscriptfont\bffam=\fivebf
  \def\tt{\fam\ttfam\eighttt}%
  \textfont\ttfam=\eighttt
  \tt \ttglue=.5em plus.25em minus.15em
  \normalbaselineskip=9pt
  \def\MF{{\manual opqr}\-{\manual stuq}}%
  \let\big=\eightbig
  \setbox\strutbox=\hbox{\vrule height7pt depth2pt width\z@}%
  \normalbaselines\rm}
\def\eightbig#1{{\hbox{$\textfont0=\ninerm\textfont2=\ninesy
  \left#1\vbox to6.5pt{}\right.\n@space$}}}

\catcode`@=13 
                                               
\def\({\left(}              %
\def\){\right)}             %

\def\A{{\cal A}}

\def\V{{\cal V}}
\def\L{\Lambda}

\cl{\bffg Transfinite dimensions}
\mp\sp
\cl{\bfff  Gerald Kuba}
\bp\mp

Let $\,|M|\,$ denote the cardinal number (the {\it size}) of the set $\,M\,$,
e.g. 
$\;|\N|=\aleph_0\,$, $\;|\R|=2^{\aleph_0}\,$.
\sp
\vbox{\eightpoint
Recall that $\;|A|^{|B|}\;$ is the
cardinal number of the family of all functions from $\,B\,$ to $\,A\,$.
Note further that if $\,M\,$ is an infinite set then $\;|M|\;$
is the size of the family of all {\it finite} subsets of $\,M\,$
and  $\;|M|^{\aleph_0}\;$
is the size of the family of all {\it countable} subsets of $\,M\,$
while, of course,  $\;2^{|M|}\;$
is the size of the family of {\it all} subsets of $\,M\,$
and $\;|M|< 2^{|M|}\,$.
Consequently, $\;|M|\leq |M|^{\aleph_0}\leq 2^{|M|}\,$.}

\mp
For a (real or complex)  Hilbert space $\,{\cal H}\,$
let $\,B\,$ be a linear basis
of $\,{\cal H}\,$ and
$\,S\,$ be an orthonormal basis
of $\,{\cal H}\,$.
Of course we have $\;|B|\geq|S|\;$ and $\;|B|=|S|\;$ when $\,B\,$
or $\,S\,$ is finite.
Further
it is well-known that $\;|B|\geq 2^{\aleph_0}\;$ when
$\,S\,$ is infinite [4].
It is also known that $\;|B|^{\aleph_0}=|B|\;$
when $\,B\,$ is infinite [6].
From the main theorem in [3] it follows that 
either $\;|B|=|S|\;$ or $\;|B|=2^{|S|}\,$.
Unfortunately, this conclusion 
is {\it not provable} in standard set theory!
Actually, the "proof$\,$" of the statement in [3] Theorem 1 uses a very strong
cardinal hypothesis (namely 
the {\it Generalized Continuum Hypothesis})
which is not mentioned in the statement.
Even worse, [3] Theorem 1  is {\it false}
under the irrefutable assumption $\;2^{\aleph_0}>\aleph_1\;$ because 
then one immediately obtains the wrong conclusion $\;|B|=|S|\;$ when
$\;|S|=\aleph_1\,$. 
\mp
Fortunately, one need not assume an unprovable hypothesis
in order to derive a natural and simple
relation between the {\it linear dimension} $\;\beta=|B|\;$
and the  {\it orthonormal dimension} $\;\sigma=|S|\;$
of an arbitrary Hilbert space $\,{\cal H}\,$.
This relation between these two fundamental concepts of {\it dimension},
which can hardly be found in the literature, reads as follows. 
\mp
{\bf Theorem 1.} {\it For every infinite-dimensional Hilbert space 
$\,{\cal H}\,$ we have 
$\;\;\beta\,=\,|{\cal H}|\,=\,\sigma^{\aleph_0}\;$.}
\mp
{\it Remark.} Concerning {\it Banach spaces} as considered in [3]
our Theorem 1 remains true when $\,{\cal H}\,$ is a 
(real or complex) Banach space with linear dimension $\,\beta\,$ and a
{\it Schauder basis} of size $\,\sigma\,$. In fact,
{\it this} has been proved but unfortunately not stated in [3].
\mp
As a consequence of Theorem 1, 
two non-isomorphic Hilbert spaces $\,{\cal H}_1\,$
and $\,{\cal H}_2\,$ can be 
isomorphic as pure vector spaces. 
For example let $\;{\cal H}_1=\ell^2=\ell^2(\N)\;$
and $\;{\cal H}_2=\ell^2(\R)\,$.
(As usual, $\,\ell^2(\L)\,$ consists of all mappings $\;x:\,\L\to \C\;$ 
such that $\,x(\lambda)\not=0\,$
for only countably many $\,\lambda\in\L\,$ and
$\;\sum_{\lambda\in\L}|x(\lambda)|^2\,<\,\infty\,$.)
Then we have $\;\beta_1=\beta_2=2^{\aleph_0}\;$ and
$\;\sigma_1={\aleph_0}<2^{\aleph_0}=\sigma_2\,$.
Further examples are all pairs
$\;{\cal H}_1=\ell^2(\Lambda_1)\,,\;{\cal H}_2=\ell^2(\Lambda_2)\;$
where $\;|\Lambda_1|=\kappa\,,\;|\Lambda_2|=\kappa^{\aleph_0}\;$
and $\;\kappa=\aleph_{\alpha+\omega}\;$
with an arbitrary ordinal $\,\alpha\,$.
(Then $\;\sigma_1=\kappa<\kappa^{\aleph_0}=\sigma_2=\beta_1=\beta_2\,$.)
\bp
If $\,\V\,$ is a vector space
over a field $\,\K\,$ then let $\;\hbox{dim}\,\V\;$ 
denote its ordinary {\it dimension}.
(In particular, $\;\hbox{dim}\,{\cal H}\;$ is 
the linear dimension $\,\beta\,$ of the Hilbert space $\,{\cal H}\,$.)
For two vector spaces $\;{\cal V}_1,{\cal V}_2\;$ over one field
we write $\;{\cal V}_1\hookrightarrow{\cal V}_2\;$
if and only if 
$\,{\cal V}_2\,$  contains an  isomorphic copy of $\,{\cal V}_1\,$.
For every field $\,\K\,$ and every set $\,I\,$ the
set $\,\K^I\,$  of all functions from $\,I\,$ to $\,\K\,$
becomes a vector space over $\,\K\,$ when the algebraic structure is defined 
in the canonical way.
If 
$\,I\,$ is a {\it finite} set then, of course, $\;\hbox{dim}\,\K^I\,=\,|I|\;$
for every field $\,\K\,$. But $\;\;\hbox{\rm dim}\,\K^I\,>\,|I|\;\;$
when $\,I\,$ is infinite. 
Moreover, the dimension of the space $\,\K^I\,$ can be
determined precisely.
\mp
{\bf Theorem 2.} {\it If  $\,{\cal V}\,$ is an infinite vector space
over an arbitrary field 
$\,\K\,$ such that $\;\K^{\N}\hookrightarrow {\cal V}\;\,$
or $\,\;|{\cal V}|>|\K|\;\;$ then
$\;\;\hbox{\rm dim}\,{\cal V}\,=\,|{\cal V}|\;$.}
\mp
{\bf Corollary 1 [}$\,${\sl The  {\smc Erd\"os}-{\smc Kaplansky} 
Theorem}$\,${\bf ].}                                           

{\it For every field $\,\K\,$ and every infinite set $\,I\,$,
$\;\;\hbox{\rm dim}\,\K^{I}\,=\,|\K|^{|I|}\;$.}
\eject
\mp
{\it Remark.} A
proof of the {\smc Erd\"os}-{\smc Kaplansky} 
Theorem is sketched in [2, Ch.II {\S}7 Ex.3]
but unfortunately this is done in an  unnecessarily complicated way.
\mp
{\bf Corollary 2.} 
{\it For every field $\,\K\,$ and every infinite set $\,I\,$, 
$\;\;\hbox{\rm dim}\,\K^{I}\,=\,2^{|I|}\;\;$ 
provided that
$\;|\K|\leq 2^{|I|}\,$.
In particular, $\;\hbox{\rm dim}\,\R^{\N}=2^{\aleph_0}\;$ 
{\it and} $\;\hbox{\rm dim}\,\R^{\R}=2^{2^{\aleph_0}}$.} 
\mp\sp
{\it Remark.} The inverse of Theorem 1 is false:
For the polynomial ring $\,{\cal V}=\Q[X]\,$ we have
$\,\;\hbox{dim}\,{\cal V}=|{\cal V}|=|\Q|=\aleph_0\,\;$ and (hence
by Corollary 1) $\,\Q^{\N}\,$ cannot be embedded in $\,{\cal V}\,$. 
\bp
First we prove Theorem 2. In doing so the following
proposition is essential.
\mp
{\bf Proposition 1.} {\it For every field $\,\K\,$ the set
\sp
\cl{$\{\,(a^n)_{n\in\N}\;\;|\;\;0\not=a\in\K\,\}\;,$}
\sp
which is obviously equipollent to $\;\K\setminus\{0\}\;$ 
and hence equipollent to $\,\K\,$ when $\,\K\,$
is infinite,
is a linearly independent subset of the vector space 
$\,\K^{\N}\,$.}
\mp

{\it Proof.} Let $\;a_1,a_2,...,a_m\;$ be distinct elements of $\,\K\,$.
We are done by verifying that for arbitrary 
$\;\lambda_1,\lambda_2,...,\lambda_m\in\K\;$ 
$$\sum_{k=1}^m\lambda_k\cdot(a_k^n)_{n\in\N}$$
is the zero vector $\;(0,0,0,...\,)\;$ in the space $\,\K^{\N}\,$
only in the trivial case $\;\lambda_1=\lambda_2=\cdots=\lambda_m=0\,$.
In other words we have to verify that 
$$\forall\,n\in\N\;:\quad \sum_{k=1}^m\lambda_k\cdot a_k^n\;=\;0$$
implies $\;\lambda_1=\lambda_2=\cdots=\lambda_m=0\,$.
Now, $\;\lambda_1=\lambda_2=\cdots=\lambda_m=0\;$
is already a consequence of the weaker assumption 
$$\forall\,n\in\{0,1,...,m\!-\!1\}\;:
\quad \sum_{k=1}^m\lambda_k\cdot a_k^n\;=\;0\;$$
because this system of equations has the matrix
$$\(\matrix{ 1 & 1 & \cdot& \cdot & \cdot & \cdot & 1 \cr
a_1 & a_2 & \cdot & \cdot & \cdot& \cdot  & a_m \cr
a_1^2 & a_2^2 & \cdot & \cdot & \cdot & \cdot & a_m^2 \cr
\cdot & \cdot & \cdot & \cdot & \cdot & \cdot & \cdot \cr
\cdot & \cdot & \cdot & \cdot & \cdot & \cdot & \cdot \cr
\cdot & \cdot & \cdot & \cdot & \cdot & \cdot & \cdot \cr
a_1^{m-1} & a_2^{m-1} & \cdot & \cdot & \cdot & \cdot & a_m^{m-1} \cr}\)$$
and, naturally, the determinant of this matrix equals  
$$\prod\limits_{1\leq i<j\leq m}\!\!\!(a_j-a_i)\;\;\not=\;\;0\;.$$ 
\bp        
{\it Proof of Theorem 2.} If $\;|{\cal V}|>|\K|\;$ then we must have
$\;\hbox{\rm dim}\,{\cal V}\,=\,|{\cal V}|\;$ because
the assumption 
$\;\hbox{\rm dim}\,{\cal V}\,<\,|{\cal V}|\;$
would immediately lead to a contradiction 
since a vector space is always equal to the union of
all subspaces which are generated by
finite subsets of a basis of the vector space.
(Note that $\,{\cal V}\,$ is assumed to be an infinite set
whence $\,{\cal V}\,$ is infinite-dimensional if $\,\K\,$
is a finite field.)
Now assume $\;|{\cal V}|=|\K|\,$. Then $\,\K\,$ is infinite and
by Proposition 1 we have 
$\;\hbox{\rm dim}\,\K^{\N}\,\geq\,|\K|\;$ 
and thus from 
$\;\K^{\N}\hookrightarrow {\cal V}\;$
we obtain $\;\hbox{\rm dim}\,{\cal V}\,\geq\,|\K|\;$
which means $\;\hbox{\rm dim}\,{\cal V}\,\geq\,|{\cal V}|\,$.
Therefore, $\;\hbox{\rm dim}\,{\cal V}\,=\,|{\cal V}|\,$.
\bp
It is worth mentioning that Corollary 2 can be proved directly
without using ideas from the proof of Theorem 2. 
\mp
{\it Proof of Corollary 2.}  
For $\;A\subset I\;$ let $\;{\bf 1}_A\;$ be the characteristic function
of the set $\,A\,$. So $\;{\bf 1}_A(x)\in\{0,1\}\subset\K\;$
for all $\;x\in I\;$ where 
$\;{\bf 1}_A(x)=1\;$ when $\;x\in A\;$ and
$\;{\bf 1}_A(x)=0\;$ when $\;x\not\in A\,$.
Let $\,{\cal F}\,$
be a family of subsets of $\,I\,$ 
with $\;|{\cal F}|=2^{|I|}\;$ such that 
$\;J\,\not\subset\,I_1\cup\cdots\cup I_n\;$
whenever $\,J\in{\cal F}\,$ and $\;I_1,...,I_n\,\in\,{\cal F}\setminus\{J\}\;$
for arbitrary $\,n\in\N\,$. (Such a family $\,{\cal F}\,$ exists by 
[5] Lemma 7.7.) Then $\,n\,$ vectors in the space $\,\K^I\,$
taken from the family $\;{\cal G}\,:=\,\{\,{\bf 1}_F\;|\;F\in{\cal F}\,\}\;$
must be linearly independent for arbitrary $\,n\in\N\,$.
This is true simply because 
the standard basis of the vector space $\,\K^n\,$
consists of $\,n\,$ linearly independent vectors.
Hence $\;\hbox{\rm dim}\,\K^I\,\geq\,|{\cal G}|=|{\cal F}|=2^{|I|}\,$. 
This is enough since $\,|\K^{I}|\leq (2^{|I|})^{|I|}=2^{|I|}\,$. 
\bp
{\it Remark.} Of course, $\;\K^{\N}\hookrightarrow{\cal V}\;$ if and only if
$\;\hbox{dim}\,{\cal V}\,\geq\,\hbox{dim}\,\K^{\N}\,=\,
|\K|^{\aleph_0}$. But consider the space                
$\;{\cal V}\,=\,\K^{I}\;$
where $\,\K\,$ is a field of size $\,\aleph_3\,$
and $\;|I|=\aleph_1\,$.
It is obvious that $\;\K^{\N}\hookrightarrow{\cal V}\;$
but it seems not possible to verify
$\;\hbox{dim}\,{\cal V}\geq|\K|^{\aleph_0}\;$ directly
without using $\;\K^{\N}\hookrightarrow{\cal V}\;$
because the standard system of linearly independent vectors in $\,{\cal V}\,$
has size $\;|I|=\aleph_1\;$ and a construction of 
$\,\aleph_3\,$ linearly independent vectors 
(which is possible by using $\;\K^{\N}\hookrightarrow{\cal V}\;$
and Proposition 1) is not enough since it is undecidable whether
$\,\aleph_3^{\aleph_0}\,$ 
is greater than or equal to $\,\aleph_3\,$. 
In view of the previous proof we have 
$\;\hbox{\rm dim}\,\K^I\,\geq\,2^{\aleph_1}\;$ but for our purpose this
is also not enough since $\,2^{\aleph_1}=\aleph_2^{\aleph_1}\,$ 
and one cannot rule out $\;\aleph_2^{\aleph_1}<\aleph_3^{\aleph_0}\,$. 
Besides, the second assumption $\;|{\cal V}|>|\K|\;$ of Theorem 2 is of no use
in order to compute $\;\hbox{dim}\,{\cal V}\;$ because 
$\;|{\cal V}|>|\K|\;$ means $\;\aleph_3^{\aleph_1}>\aleph_3\;$ 
which is unprovable.
\bp\bp
It remains to prove Theorem 1. To begin with we extend Theorem 2.
\mp\sp
{\bf Theorem 3.} {\it If  $\,{\cal V}\,$ is a real or complex vector space
such that the Hilbert space $\,\ell^2\,$ can be algebraically embedded
into $\,{\cal V}\,$ then $\;\;\hbox{\rm dim}\,{\cal V}\,=\,|{\cal V}|\;$.
In particular, $\;\hbox{\rm dim}\,{\ell^2}\,=\,2^{\aleph_0}\,$.}
\mp\sp
The proof of Theorem 3 is a simple adaption of the proof of Theorem 2
in view of the fact that the set 
$\;\,\{\,(a^n)_{n\in\N}\;\;|\;\;0\not=|a|<1\,\}\;\,$
is a subset of $\,\ell^2\,$ of size $\,2^{\aleph_0}\,$
which (by applying Proposition 1) is linearly independent.

\mp\sp
{\it Remark.} Theorem 3 is useful to compute the dimension 
of several Banach spaces. For example, the dimension of the space of all
bounded functions from an arbitrary infinite set $\,X\,$ to $\,\R\,$ 
(or to $\,\C\,$) is equal to $\,2^{|X|}\,$.
Further it is a nice exercise to construct an embedding
of $\,\ell^2\,$ in order to show that
the dimension of the space of all
continuous functions from $\,[0,1]\,$ to $\,\R\,$ is equal to
$\,2^{\aleph_0}\,$.
\vfill\eject
\bp
With the help of Theorem 3 we immediately obtain Theorem 1.
Let $\,{\cal H}\,$ be an infinite-dimensional Hilbert space.
Then, naturally, 
$\;\ell^2\hookrightarrow {\cal H}\;$
and hence $\;\beta\,=\,\hbox{\rm dim}\,{\cal H}\,=\,|{\cal H}|\;$.
It remains to verify $\;|{\cal H}|=\sigma^{\aleph_0}\;$
where $\,\sigma\,$ is the orthonormal dimension of 
$\,{\cal H}\,$. Since $\,{\cal H}\,$ is {\it norm-isomorphic}
to $\;\ell^2(\Lambda)\;$ with $\;|\Lambda|=\sigma\;$
the proof of Theorem 1 is finished by showing 
the following proposition.
\mp\sp
{\bf Proposition 2.} {\it 
$\;\;|\ell^2(\Lambda)|=|\Lambda|^{\aleph_0}\;$
for every infinite index set $\,\Lambda\,$.}
\mp
{\it Proof.} Let $\,\K\,$ be either the field $\,\R\,$ or the field
$\,\C\,$. By definition,
\sp
\cl{$\;\ell^2(\L)\,=\,\big\{\,x\in\A(\L,\K)\,\,\big|\,\,
\sum\limits_{\lambda\in\L}x(\lambda)^2\,<\,\infty\,\,\big\}\;$}
\sp
where $\,\A(\L,\K)\,$ is the family of all
functions \hbox{$\;x:\,\L\to \K\;$} such that $\,x(\lambda)\not=0\,$
for at most countably many $\,\lambda\in\L\,$.
\sp
Elementary transfinite arithmetics yields
$\;|\A(\Lambda,\K)|=|\Lambda|^{\aleph_0}|\K|^{\aleph_0}=
|\L|^{\aleph_0}\,$. Consequently, 
$\;|\ell^2(\L)|\leq |\L|^{\aleph_0}\,$.
It remains to verify $\;|\ell^2(\L)|\geq |\L|^{\aleph_0}\,$.
\sp
Let $\,{\cal I}(\L)\,$ be the family of all injective functions
from $\,\N\,$ to $\,\L\,$. Naturally, $\;|{\cal I}(\L)|=|\L|^{\aleph_0}\,$.
Obviously, the space $\,\ell^2(\L)\,$ contains the set
\sp
\cl{${\cal F}(\L)\;:=\;\big\{\,x\in\R^\Lambda\,\,\big|\,\,
\exists\,\varphi\in{\cal I}(\L)\;:\,\,
x(\L\!\setminus\!\varphi(\N))=\{0\}\,\,\land\,\,
\big(x(\varphi(n))\big)_{n\in \N}\,\in\,
\prod\limits_{n=1}^\infty[0,{1\over n}]\,\big\}\,$.}

From $\;|{\cal I}(\L)|=|\L|^{\aleph_0}\;$                   
and $\;\big|\!\prod\limits_{n=1}^\infty [0,{1\over n}]\,\big|=2^{\aleph_0}\;$ 
we derive $\;|{\cal F}(\L)|=|\L|^{\aleph_0}\,$, {\it q.e.d.}
\bp
{\it Remark.} Due to Anderson's theorem [1], Hilbert space
$\,\ell^2=\ell^2(\N)\,$ is homeomorphic to the product space 
$\,\R^{\N}\,$. And by Theorem 1 and Corollary 2 also the vector spaces 
$\,\ell^2\,$ and $\,\R^{\N}\,$ are isomorphic. 
However, Proposition 2 demonstrates that neither the first nor the second 
statement can be generalized from $\,\L=\N\,$ to arbitrary index sets $\,\L\,$. 
(If $\,|\Lambda|=2^\kappa\,$ for some transfinite 
cardinal $\,\kappa\,$
then $\;|\ell^2(\Lambda)|=|\Lambda|^{\aleph_0}=|\L|<2^{|\L|}=|\R^\L|\;$
and hence $\,\ell^2(\L)\,$ cannot be homeomorphic or 
algebraically isomorphic to $\,\R^{\L}\,$.)
\bp\bp\mp
{\bff References}
\bp
[1] Anderson, R.D.: {\it Hilbert space is homeomorphic to the 
countable infinite product of}

\rl{{\it lines.} Bull.~Amer.~Math.~Soc.~{\bf 72}, 
515-519 (1966).}
\sp
[2] Bourbaki: Elements of Mathematics, Algebra I (2nd printing).
Springer 1989.
\mp
[3] Evans, J.W., and Tapia, R.A.: {\it Hamel versus Schauder dimension.}

\rl{Am.~Math.~Monthly~{\bf 77}, No.4, 385-388 (1970).}
\sp
[4] Halmos, P.R.: A Hilbert Space Problem Book. Springer 1974.
\mp
[5] Jech, Th.: Set Theory. Third Millennium Edition. Springer 2002.
\mp
[6] Kruse, A.H.: {\it Badly incomplete normed linear spaces.} 
Math.~Z.~{\bf 83}, 314-320 (1964).

\bp\bp
{\sl Author's address:} Institute of Mathematics. 

University of Natural Resources and Life Sciences, Vienna, Austria. 

{\sl E-mail:} {\tt gerald.kuba(at)boku.ac.at}
\end
\bp\bp
{\bff Author's address} 
\mp
Institute of Mathematics
\sp
University of Natural Resources and Life Sciences 
\sp
Vienna, Austria. 
\mp
{\sl E-mail:} {\tt gerald.kuba(at)boku.ac.at}
\end

{\sl Author's address:} Institute of Mathematics. 

University of Natural Resources and Life Sciences, Vienna, Austria. 

{\sl E-mail:} {\tt gerald.kuba@boku.ac.at}